\documentclass[12pt,reqno]{article}
\usepackage{amssymb}
\usepackage{amsthm}
\usepackage{amsmath}
\newtheorem{theorem}{Theorem}[section]
\newtheorem{corollary}[theorem]{Corollary}

\newtheorem{example}[theorem]{Example}

\newtheorem{remark}[theorem]{Remark}

\def\be{\begin{equation}}
\def\ee{\end{equation}}
\def\bea{\begin{eqnarray}}
\def\eea{\end{eqnarray}}
\begin{document}
\begin{center} \Large{\bf Linearization of third-order ordinary differential
equations $u'''=f(x,u,u',u'')$ via point transformations}
\end{center}
\medskip
\begin{center}
Ahmad Y. Al-Dweik$^*$, M. T. Mustafa$^{**}$,
F. M. Mahomed$^{***}$ and R. S. Alassar$^*$\\

{$^*$Department of Mathematics \& Statistics, King Fahd University
of Petroleum and Minerals, Dhahran 31261, Saudi Arabia}\\
{$^{**}$Department of Mathematics, Statistics and Physics, Qatar
University, Doha, 2713, State of Qatar}\\
{$^{***}$School of Computer Science and Applied Mathematics, DST-NRF
Centre of Excellence in Mathematical and Statistical Sciences,
University of the Witwatersrand, Johannesburg, Wits 2050,  South Africa\\
}

aydweik@kfupm.edu.sa, tahir.mustafa@qu.edu.qa,
Fazal.Mahomed@wits.ac.za and alassar@kfupm.edu.sa.
\end{center}
\begin{abstract}
The linearization problem by use of the Cartan equivalence method
for scalar third-order ODEs via point transformations was solved
partially in \cite{DweikTahirFazal2,DweikTahirFazal1}. In order to solve
this problem completely, the Cartan equivalence method is applied
to provide an invariant characterization of the linearizable
third-order ordinary differential equation $u'''=f(x,u,u',u'')$
which admits a four-dimensional point symmetry Lie algebra. The
invariant characterization is given in terms of the function $f$
in a compact form. A simple procedure to construct the equivalent
canonical form by use of an obtained invariant is also presented.
The method provides auxiliary functions which can be utilized to
efficiently  determine the point transformation that does the
reduction to the equivalent canonical form. Furthermore, illustrations to the main theorem and
applications are given.
\end{abstract}
\section{Introduction}
The Lie algebraic criteria for linearization  for scalar $n$th-order ($n>2$) ordinary differential equations (ODEs)
by means of point transformations were uncovered
in \cite{mah4}. The canonical forms for scalar third-order ODEs are listed in
\cite{Mahomed1996}. Three canonical forms occur for scalar linear
third-order ODEs. The maximal Lie algebra case for such ODEs is of
dimension seven and corresponds to the simplest equation $u'''=0$.

The Laguerre-Forsyth (see \cite{Mahomed2007}) canonical form for scalar
linear third-order ODEs is
\begin{equation}\label{eee3}
u'''+a(x)u=0.
\end{equation} If $a\not  \equiv 0$, then (\ref{eee3}) has a five- or four-dimensional symmetry Lie
algebra.

Chern \cite{che} pioneered the use of the Cartan equivalence
method to solve the linearization problem of scalar third-order
ODEs via contact transformations. He determined conditions for
equivalence to (\ref{eee3}) in the cases $a  \equiv0$ and $a
\equiv1$. Neut and Petitot \cite{neu} in contemporary times then
obtained criteria for equivalence to (\ref{eee3}) by contact
transformations for arbitrary $a(x)$. Grebot \cite{gre}
investigated linearization of third-order ODEs via fibre
preserving invertible transformations. In recent work, Ibragimov
and Meleshko \cite{ibr} studied the linearization problem by
utilizing a direct Lie approach for such third-order ODEs by both
point and contact transformations. These authors also presented
conditions on the linearizing maps. In the works
\cite{DweikTahirFazal2,DweikTahirFazal1}, the authors very recently applied the
Cartan equivalence method  to deduce an invariant characterization
of the scalar third-order ODE $u'''=f(x,u,u',u'')$ which possesses
both five- and seven-dimensional point symmetry Lie algebras.
Moreover, we provided auxiliary functions which can be utilised to
efficiently   obtain the point transformation.

Therefore, in \cite{DweikTahirFazal2,DweikTahirFazal1}, the authors have partially
solved the linearization problem for third-order ODEs $u'''=f(x,u,u',u'')$ via invertible point
transformations
\begin{equation}\label{cccc}
\bar{x}=\phi \left( x,u \right),~\bar{u} =\psi \left( x,u  \right),~~~\phi_x\psi_u   -  \phi_u \psi_x\neq0\\
\end{equation}
by the Cartan equivalence method.
In order to complete the previous studies, we investigate herein the
linearization problem  by the Cartan approach for
scalar third-order ODEs, via invertible point transformation, which admit a four dimension point symmetry
algebra.

It is prudent to mention here that the invariant characterization of
$u'''=f(x,u,u',u'')$ that admits both five and seven  point
symmetry algebra was found in terms of the function $f$ in the
following theorems. We denote $u', u''$ by $p, q$ respectively, in
the following and hereafter.
\begin{theorem}\label{TH1.1}\cite{DweikTahirFazal2}.
The necessary and sufficient conditions for equivalence of a
scalar third-order ODE $u''' = f(x,u,u',u'')$ via {\it point
transformation} (\ref{cccc}) to its canonical
form $\bar{u}'''=0$, with seven point symmetries, are the identical vanishing of the
relative invariants
\begin{equation}\label{bbbb40}
\begin{array}{l}
  I_1  = f_{q,q,q}  \hfill \\
  I_2  = f_{q,q} ^2  + 6\,f_{p,q,q}  \hfill \\
  I_3  =  {4\,f_q ^3  + 18\,f_q \left( {f_p  - D_x f_q } \right) + 9\,D_x^2 f_q    - 27\,D_x f_p+ 54\,f_u  }\\
  I_4 =f_{q,q} \left( {f_q ^2+9\,f_p  - 3\,D_x f_q}\right) - 9\,f_{p,p}  + 18\,f_{u,q}  - 6\,f_q f_{p,q}.\\
\end{array}
\end{equation}
Given that the the system of relative invariants (\ref{bbbb40}) is
zero, the linearizing point transformation (\ref{cccc}) is defined
by
\begin{equation}\label{bbbb42}
\begin{array}{l}
D_x\phi=\frac{a_1}{a_3},\\
\phi_x\psi_u   - \phi_u \psi_x=\frac{a^2_1}{a_3},\\
\end{array}
\end{equation}
where $a_1(x,u,p), a_2(x,u,p,q),a_3(x,u,p)$ are auxiliary
functions given by
\begin{equation}\label{bbbb43}
\begin{array}{l}
D_xa_3=-\frac{1}{3}f_q a_3,\\
D_xa_2=\frac{1}{2}\frac{1}{a_3}a^2_2-\frac{1}{18}a_3\left( 2f_q ^2 +9\,f_p   - 3\,D_x f_q \right),\\
D_x a_1=\frac{a_2}{a_3}a_1,\\
\left(\frac{a_1}{a_3}\right)_{p,p}=0,\\
\left(\frac{a^2_1}{a_3}\right)_{p}=0.\\
\end{array}
\end{equation}
\end{theorem}
\begin{remark}\label{r1}\cite{DweikTahirFazal2}.
The system (\ref{bbbb43}) provides
\begin{equation}\label{b44}
\begin{array}{l}
a_2(x,u,p,q)=-\frac{1}{6}a_3 f_{q,q}~q+A(x,u,p),
\end{array}
\end{equation}
for some function $A(x,u,p)$.
\end{remark}
\begin{theorem}\label{TH1.2}\cite{DweikTahirFazal1}.
The necessary and sufficient conditions for equivalence of a
third-order equation $ u''' = f(x,u,u',u'')$ via {\it point transformation} (\ref{cccc}) to the canonical form
$\bar{u}'''=s~\bar{u}'+\bar{u},~s=\,{\rm constant}$, with five
point symmetries,  are
the identically vanishing of the relative invariants
\begin{equation}\label{bsd40}
\begin{array}{l}
  I_1  = f_{q,q,q}  \hfill \\
  I_2  = f_{q,q} ^2  + 6\,f_{p,q,q}  \hfill \\
  I_4  = J_q \,\,\,\,\,\,\,\,\,\,\,\,\,\,\,\,\,\,\,\,\,\,\,\,\,\,\,\,\,\,\,\,\,\,\,\,\,\,\,\,\,\,\,\,\,\,\,\,\,\,\,\,\,\,\,\,\,\,\,\,\,\,\,\,\,\,\,\,\,\,\,\,\,\,\,\,\,\,\,\,\,\,\,\,\,\,\,\,\,\,\,\,\,\,\,\,\,\,\,\,\, \hfill \\
  I_5  = f_{q,q} \,J - 6\,J_p \,\,\,\,\,\,\,\,\,\,\,\,\,\,\,\,\,\,\,\,\,\,\,\,\,\,\,\,\,\,\, \hfill \\
  I_6  = J_u  - D_x J_p  \hfill \\
  I_7  = f_{q,q} \left( {9\,f_p  + f_q ^2  - 3\,D_x f_q } \right) - 9\,f_{p,p}  + 18\,f_{u,q}  - 6\,f_q f_{p,q}  \hfill \\
  I_9  =K_q\\
  I_{10}=K_p\\
  I_{11}=K_u\\
  I_{12}=K_x,\\
\end{array}
\end{equation}
where
\begin{equation}\label{bsd41}
\begin{array}{l}
I_3  =J^3= \frac{1}{54} \left( {4\,f_q ^3  + 18\,f_q \left( {f_p  - D_x f_q } \right) + 9\,D_x^2 f_q    - 27\,D_x f_p+ 54\,f_u  }\right)  \ne 0,\\
I_8 = \frac{1}{3} \left( \left( {f_q ^2  + 3\,f_p  - 3\,D_x f_q \,} \right)\,J^2 \, + 6\,J\,D_x^2 J - 9\,\left( {D_x J} \right)^2\right), \\
K=\frac{I_8}{J^4}.\\
\end{array}
\end{equation}
Provided that the the system of relative invariants (\ref{bsd40}) is
zero, the linearizing point transformation (\ref{cccc}) is defined
by
\begin{equation}\label{bsd42}
\begin{array}{l}
D_x \phi=J, \\
\phi_x\psi_u- \phi_u \psi_x=J a_1,\\
\end{array}
\end{equation}
where the auxiliary function $a_1(x,u,p)$ satisfies
\begin{equation}\label{bsd43}
\begin{array}{l}
\frac{D_x a_1}{a_1}=\frac{1}{3}\left(\frac{3 D_x J-J f_q}{J}\right).\\
\end{array}
\end{equation}
Finally, the constant $s$ of the resulting canonical form is given
by the equation $s=K$.
\end{theorem}
\begin{corollary}
The $\mathbf{linearizable}$ third-order ODEs with four point
symmetries have the Laguerre-Forsyth canonical form
\begin{equation}\label{a1}
\begin{array}{lll}
 u'''= a^3( x)  u,& a( x)\ne0,& \left(\frac{2aa''-3{a'}^2}{a^4}\right)_x\ne0.\\
\end{array}
\end{equation}
\end{corollary}
\proof It is well-known that linear third-order ODEs can be
transformed to the Laguerre-Forsyth canonical form $u'''= a^3( x)
u$ via the Laguerre transformations.

By using Theorem \ref{TH1.1}, the Laguerre-Forsyth canonical form
$u'''= a^3( x) u$ admits seven point symmetries if and only if
$a(x)=0$.

Also by utilizing Theorem \ref{TH1.2}, the Laguerre-Forsyth canonical form
$u'''= a^3( x) u$ admits five point symmetries if and only if
$a(x)\ne0$ but $\left(\frac{2aa''-3{a'}^2}{a^4}\right)_x=0$.

Finally, it is well-known that the linear third-order ODEs admit
seven, five or four point symmetries.
\endproof
In this paper, we consider the canonical form (\ref{a1}) and
Yumaguzhin's canonical  form \cite{Yum} for
$\mathbf{linearizable}$ third-order ODEs with four point
symmetries. In section 2,  the Cartan equivalence method is
applied to obtain an invariant characterization of the scalar
third-order ODE $u'''=f(x,u,u',u'')$ with four point symmetries.
Moreover, an efficient algorithm  is given to determine the point
transformation that does the reduction to the equivalent canonical
form (\ref{a1}). Section 3 provides illustrations of the result of
section 2. In section 4, we give an efficient algorithm to
determine the point transformation that does the reduction to
Yumaguzhin's canonical form. Section 5 is devoted to the
illustration of the result of section 4. In section 6, we apply
our main results to solve the point symmetry group classification
problem of third-order linear ODE.
\section{Linearizable third-order ODEs with four point symmetries}
Let $x:=(x,u,p=u',q=u'')\in \mathbb{R}^4$ be local coordinates of
$J^2$, the space of the second order jets. In local coordinates,
the equivalence of
\begin{equation}\label{b0}
\begin{array}{cc}
u'''=f(x,u,u',u''), & \bar{u}'''=\bar{f}(\bar{x},\bar{u},\bar{u}',\bar{u}''),\\
\end{array}
\end{equation}
under a point transformation (\ref{cccc}) is expressed as the
local equivalence problem for the $G$-structure
\begin{equation}\label{b2}
\Phi^*\left(%
\begin{array}{c}
  \bar{\omega}^1 \\
  \bar{\omega}^2 \\
  \bar{\omega}^3 \\
  \bar{\omega}^4 \\
\end{array}%
\right)=\left(%
\begin{array}{cccc}
  a_1 & 0 & 0 & 0 \\
  a_2 & a_3 & 0 & 0 \\
  a_4 & a_5 & a_6 & 0 \\
  a_7 & 0 & 0 & a_8 \\
\end{array}%
\right) \left(%
\begin{array}{c}
  \omega^1 \\
  \omega^2 \\
  \omega^3 \\
  \omega^4 \\
\end{array}%
\right),
\end{equation}
where
\begin{equation}\label{b1}
\begin{array}{llll}
\bar{\omega}^1=d\bar{u}-\bar{p} d\bar{x}, & \bar{\omega}^2=d\bar{p}-\bar{q} d \bar{x}, & \bar{\omega}^3=d\bar{q}-\bar{f} d \bar{x}, & \bar{\omega}^4= d \bar{x},\\
\omega^1=du-p d x, & \omega^2=dp-q d x, & \omega^3=dq-f d x, & \omega^4= d x.\\
\end{array}
\end{equation}
One can evaluate the functions $a_i=a_i(x,u,p,q), i=1\dots 8,$
explicitly. For instance, $a_1=\frac{\phi_x\psi_u   - \phi_u
\psi_x}{D_x \phi}$.

Now, we define $\theta$ to be the lifted coframe with an
eight-dimensional group $G$ as
\begin{equation}\label{b21}
\left(%
\begin{array}{c}
  \theta^1 \\
  \theta^2 \\
  \theta^3 \\
  \theta^4 \\
\end{array}%
\right)=\left(%
\begin{array}{cccc}
  a_1 & 0 & 0 & 0 \\
  a_2 & a_3 & 0 & 0 \\
  a_4 & a_5 & a_6 & 0 \\
  a_7 & 0 & 0 & a_8 \\
\end{array}%
\right) \left(%
\begin{array}{c}
  \omega^1 \\
  \omega^2 \\
  \omega^3 \\
  \omega^4 \\
\end{array}%
\right).
\end{equation}
Cartan's method, when applied to this equivalence problem, leads
to an ${e}$-structure, which is invariantly associated to the
given equation.

We note that the canonical form (\ref{a1}) of the linearizable
third-order ODEs with four point symmetries belongs to the branch
$I_1=I_2=0$, $I_3\ne0$, $I_4=I_5=I_6=0$, where
\begin{equation}\label{xxx}
\begin{array}{ll}
I_1  &=f_{q,q,q}\\
I_2  &=f_{q,q} ^2  + 6\,f_{p,q,q}\\
I_3  &=J^3= \frac{1}{54} \left( {4\,f_q ^3  + 18\,f_q \left( {f_p  - D_x f_q } \right) + 9\,D_x^2 f_q  - 27\,D_x f_p + 54\,f_u }\right) \hfill \\
I_4  &=J_q\\
I_5  &= f_{q,q} \,J - 6\,J_p \\
I_6  &= J_u  - D_x J_p \\
\end{array}
\end{equation}
are system of relative invariants derived in the work
\cite{DweikTahirFazal1}.

As it was shown in the work \cite{DweikTahirFazal1}, the {\it
fifth-order} normalizations reduces the lifted coframe (\ref{b21})
to an adapted coframe with the {\it one-dimensional group} $G_5$
\begin{equation}\label{b223}
\left(%
\begin{array}{c}
  \theta^1 \\
  \theta^2 \\
  \theta^3 \\
  \theta^4 \\
\end{array}%
\right)=\left(%
\begin{array}{cccc}
  a_1 & 0 & 0 & 0 \\
  \frac{1}{3}\frac{a_1}{J^2}s_4 & \frac{a_1}{J} & 0 & 0 \\
  \frac{1}{18}\frac{a_1}{J^4}s_4^2-\frac{a_1}{18J^2}s_2 & \frac{1}{3}\frac{a_1}{J^3}s_4-\frac{a_1}{3J^2} s_1 & \frac{a_1}{J^2} & 0 \\
  \frac{1}{6}Js_3 & 0 & 0 & J \\
\end{array}%
\right)  \left(%
\begin{array}{c}
  \omega^1 \\
  \omega^2 \\
  \omega^3 \\
  \omega^4 \\
\end{array}%
\right),
\end{equation}
where $s_1=f_q, s_2= 2f_q ^2 +9\,f_p- 3\,D_x f_q,~s_3=f_{q,q},~
s_4= 3 D_x J-J f_q$.

This results in the structure equation
\begin{equation}\label{b23}
d\left(%
\begin{array}{c}
  \theta^1 \\
  \theta^2 \\
  \theta^3 \\
  \theta^4 \\
\end{array}%
\right)=\left(%
\begin{array}{cccc}
  \alpha_1 & 0 & 0 & 0 \\
  0 & \alpha_1 & 0 & 0 \\
  0 & 0 & \alpha_1 & 0 \\
  0 & 0 & 0 & 0 \\
\end{array}%
\right)\wedge \left(%
\begin{array}{c}
  \theta^1 \\
  \theta^2 \\
  \theta^3 \\
  \theta^4 \\
\end{array}%
\right)+
\left(%
\begin{array}{c}
  -\theta^2 \wedge \theta^4 \\
 T^2_{14}~\theta^1 \wedge \theta^4-\theta^3 \wedge \theta^4  \\
 T^3_{12}~\theta^1 \wedge \theta^2-\theta^1 \wedge \theta^4+T^2_{14}~\theta^2 \wedge \theta^4  \\
 0 \\
\end{array}%
\right).\\
\end{equation}
where the invariants $T^2_{14}$ and  $T^3_{12}$
\begin{equation}\label{b25}
\begin{array}{ll}
T^2_{14}=-\frac{I_8}{2~J^4},\\
T^3_{12}=-\frac{I_7}{9a_1~J},\\
\end{array}
\end{equation}
are given in terms of the following relative invariants
\begin{equation}\label{b26}
\begin{array}{ll}
  I_7  = f_{q,q} \left(f_q ^2+ 9\,f_p    - 3\,D_x f_q  \right) - 9\,f_{p,p}  + 18\,f_{u,q}  - 6\,f_q f_{p,q},  \hfill \\
  I_8 = \frac{1}{3} \left( \left( {f_q ^2  + 3\,f_p  - 3\,D_x f_q \,} \right)\,J^2 \, + 6\,J\,D_x^2 J - 9\,\left( {D_x J} \right)^2\right) . \\
\end{array}
\end{equation}
It should be remarked here that the relative invariant  $I_7=0$ and
that the invariant
\begin{equation}\label{z0}
\begin{array}{ll}
 K=\frac{I_8}{J^4}=\frac{2a(x)a''(x)-3{a'(x)}^2}{a(x)^4}\ne0,\\
\end{array}
\end{equation}
for the canonical form (\ref{a1}).
\section*{Branch 1. $I_7=0$ and $I_8\ne0$}
In this branch, there are no group-dependent invariants among the
remaining unabsorbable torsion (\ref{b25}). Hence, one needs
to be careful that there are no group-dependent invariants among
coframe derivatives of the torsion (\ref{b25}) before checking the
involutivity. Dual to the lifted coframe (\ref{b223}) are the
invariant derivations
\begin{equation}\label{z1}
\begin{array}{ll}
  \frac{\partial}{\partial \theta^1}  =-\frac{s_3}{6a_1}D_x+\frac{1}{a_1}\frac{\partial}{\partial u}-\frac{s_4}{3Ja_1}\frac{\partial}{\partial p}+\frac{1}{18J^2a_1}(s_2J^2-2s_1s_4J+s^2_4)\frac{\partial}{\partial q},  \\
  \frac{\partial}{\partial \theta^2}  =\frac{J}{a_1}\frac{\partial}{\partial p}+\frac{1}{3a_1}(s_1J-s_4)\frac{\partial}{\partial q}, \\
  \frac{\partial}{\partial \theta^3}  =\frac{J^2}{a_1}\frac{\partial}{\partial q}, \\
  \frac{\partial}{\partial \theta^4}  =\frac{1}{J}D_x. \\
\end{array}
\end{equation}
It should be noted here that the invariants
\begin{equation}\label{z2}
\begin{array}{ll}
 \frac{\partial K}{\partial \theta^1}=\frac{\partial K}{\partial \theta^2}= \frac{\partial K}{\partial \theta^3}=0, \frac{\partial K}{\partial \theta^4}\ne0,\\
\end{array}
\end{equation}
hold for the canonical form (\ref{a1}). Thus, we choose the following
branch.
\section*{Branch 1.1. $K_q=K_p=f_{q,q} D_xK-6K_u=0$, $D_xK\ne0$}
In this branch, there are no group-dependent invariants among the
non-zero coframe derivative $\frac{\partial K}{\partial\theta^4}$,
so the final remaining group variable $a_1$ cannot be normalized.
Moreover, $\alpha_1$ is now uniquely defined, so the problem is
determinant. This results in the following $e$-structure on the
five-dimensional prolonged space $M^{(1)}=M \times G_5$
\begin{equation}\label{b36}
\begin{array}{llll}
  \theta^1=a_1\omega^1,\\
  \theta^2=\frac{1}{3}\frac{a_1}{J^2}s_4\omega^1+\frac{a_1}{J}\omega^2,\\
  \theta^3=\frac{1}{18}(\frac{a_1}{J^4}s_4^2-\frac{a_1}{J^2}s_2)\omega^1+ \frac{1}{3}(\frac{a_1}{J^3}s_4-\frac{a_1}{J^2} s_1)\omega^2+\frac{a_1}{J^2}\omega^3,\\
  \theta^4=\frac{1}{6}Js_3\omega^1+J \omega^4,\\
  \theta^5=\alpha_1=\frac{da_1}{a_1}+\frac{1}{36 J}(4s_3 s_4-18 {s_4}_p+J s_1 s_3-6J {s_1}_p)\omega^1+\frac{1}{6}s_3 \omega^2-\frac{1}{3}\frac{s_4}{J}\omega^4.\\
\end{array}
\end{equation}
This gives rise to the structure equations
\begin{equation}\label{b35}
\begin{array}{l}
  d\theta^1=-\theta^1 \wedge \alpha_1-\theta^2 \wedge \theta^4 \\
  d\theta^2=-\frac{K}{2}~\theta^1 \wedge \theta^4-\theta^2 \wedge \alpha_1-\theta^3 \wedge \theta^4  \\
  d\theta^3=-\theta^1 \wedge \theta^4-\frac{K}{2}~\theta^2 \wedge \theta^4-\theta^3 \wedge \alpha_1  \\
  d\theta^4=0 \\
  d\alpha_1=0. \\
\end{array}
\end{equation}
Dual to the invariant coframe (\ref{b36}) are the invariant
derivations
\begin{equation}\label{z3}
\begin{array}{ll}
  \frac{\partial}{\partial \theta^1}  =&-\frac{s_3}{6a_1}D_x+\frac{1}{a_1}\frac{\partial}{\partial u}-\frac{s_4}{3Ja_1}\frac{\partial}{\partial p}+\frac{1}{18J^2a_1}(s_2J^2-2s_1s_4J+s^2_4)\frac{\partial}{\partial q}  \\
  &+\frac{1}{36 J}\left( 18 {s_4}_p +6 J {s_1}_p-s_1 s_3 J-4 s_3 s_4\right)\frac{\partial}{\partial a_1},\\
  \frac{\partial}{\partial \theta^2}  =&\frac{J}{a_1}\frac{\partial}{\partial p}+\frac{1}{3a_1}(s_1J-s_4)\frac{\partial}{\partial q}-\frac{1}{6 }J  s_3\frac{\partial}{\partial a_1}, \\
  \frac{\partial}{\partial \theta^3}  =&\frac{J^2}{a_1}\frac{\partial}{\partial q}, \\
  \frac{\partial}{\partial \theta^4}  =&\frac{1}{J}D_x+\frac{s_4 a_1}{3 J^2}\frac{\partial}{\partial a_1}, \\
  \frac{\partial}{\partial \theta^5}  =&a_1 \frac{\partial}{\partial a_1}. \\
\end{array}
\end{equation}
By using the dual Lie brackets identities formulae (see
\cite[Equation 8.20, page 261]{Olver1995})
\begin{equation}
\left[\frac{\partial }{\partial \theta^j},\frac{\partial }{\partial \theta^k}\right]=\frac{\partial^2 }{\partial \theta^j\partial\theta^k}-\frac{\partial^2 }{\partial \theta^k\partial\theta^j}=-\sum_{i=1}^{5}T^i_{jk}\frac{\partial }{\partial \theta^i},,~1\leq j< k \leq 5,\\
\end{equation}
the invariant derivations (\ref{z3}) have the non-zero Lie
brackets
\begin{equation}\label{z4}
\begin{array}{lll}
\left[\frac{\partial }{\partial \theta^1},\frac{\partial }{\partial \theta^4}\right]=\frac{K}{2}\frac{\partial }{\partial \theta^2}+\frac{\partial }{\partial \theta^3},&\left[\frac{\partial }{\partial \theta^1},\frac{\partial }{\partial \theta^5}\right]=\frac{\partial }{\partial \theta^1},&\left[\frac{\partial }{\partial \theta^2},\frac{\partial }{\partial \theta^4}\right]=\frac{\partial }{\partial \theta^1}+\frac{K}{2}\frac{\partial }{\partial \theta^3},\\
\left[\frac{\partial }{\partial \theta^2},\frac{\partial }{\partial \theta^5}\right]=\frac{\partial }{\partial \theta^2},&\left[\frac{\partial }{\partial \theta^3},\frac{\partial }{\partial \theta^4}\right]=\frac{\partial }{\partial \theta^2},&\left[\frac{\partial }{\partial \theta^3},\frac{\partial }{\partial \theta^5}\right]=\frac{\partial }{\partial \theta^3}.\\
\end{array}
\end{equation}
In this branch, the invariant $K$ satisfies $\frac{\partial
K}{\partial \theta^1}=\frac{\partial K}{\partial \theta^2}=
\frac{\partial K}{\partial \theta^3}=\frac{\partial K}{\partial
\theta^5}=0$. Then the first-order classifying set is
\begin{equation}\label{z5}
\begin{array}{ll}
\mathcal{C}^{(1)}=\{ K,\frac{\partial K}{\partial \theta^4} \}.\\
\end{array}
\end{equation}
By utilizing the dual Lie brackets (\ref{z4}),  the $s^{th}$-order
classifying set for $s\geq2$ is
\begin{equation}\label{z6}
\begin{array}{ll}
\mathcal{C}^{(s)}=\{ K,\frac{\partial K}{\partial \theta^4},\frac{{\partial}^2 K}{\partial (\theta^4)^2},...,\frac{{\partial}^s K}{\partial (\theta^4)^s} \}.\\
\end{array}
\end{equation}
The invariants $K$ and $\frac{\partial K}{\partial\theta^4}$ are
functionally dependent, where
\begin{equation}\label{z7}
dK \wedge d\left(\frac{\partial K}{\partial\theta^4}\right)=\left(\sum_{i=1}^{5}\frac{\partial K}{\partial\theta^i} d \theta^i \right)\wedge \left(\sum_{j=1}^{5}\frac{\partial^2 K}{\partial\theta^j \partial\theta^4} d \theta^j \right)=\frac{\partial K}{\partial\theta^4} d \theta^4 \wedge \frac{{\partial}^2 K}{\partial{(\theta^4)}^2} d \theta^4=0.\\
\end{equation}
Similarly, one can show that the invariants $K$ and
$\frac{{\partial}^s K}{\partial (\theta^4)^s}$ are functionally
dependent for $s\geq2$. Therefore, if  $\rho_s$ denotes the number
of functionally independent structure invariants up to order $s$,
then $\rho_0=\rho_1=\rho_2=...=\rho_s=1$ for $s\geq2$. That is, the
stabilizing $\textit{rank}$ of  the coframe (\ref{b36}) is one and
the $\textit{order}$ of this coframe  is zero.

Therefore, we have produced an invariant coframe (\ref{b36}) with
rank one on the five-dimensional space coordinates $x,u,p, q,
a_1$. Any such differential equation admits a four-dimensional
symmetry group of point transformations (see \cite[Theorem 8.22,
page 275]{Olver1995}).

Finally, the necessary and sufficient conditions for equivalence
of a third-order equation $u''' = f(x,u,u',u'')$ via  point
transformation (\ref{cccc}) to the canonical form
$\bar{u}'''=\bar{a}^3( \bar{x})  \bar{u}$ with four point
symmetries  is that it belongs to this branch and its first-order
classifying set $\mathcal{C}^{(1)}(\theta)$ overlaps with the
first-order classifying set $\mathcal{C}^{(1)}(\bar{\theta})$ of
the canonical form $\bar{u}'''=\bar{a}^3( \bar{x})  \bar{u}$ (see
\cite[Theorem 8.19, page 271]{Olver1995}).

The overlapping of their first-order classifying sets
$\mathcal{C}^{(1)}(\theta)$ and $\mathcal{C}^{(1)}(\bar{\theta})$
gives the conditions
\begin{equation}\label{z8}
\begin{array}{l}
K=\bar{K}=\frac{2\bar{a}(\bar{x})\bar{a}''(\bar{x})-3{\bar{a}'(\bar{x})}^2}{\bar{a}(\bar{x})^4},\\
\frac{1}{\bar{J}}D_{x} K=\frac{1}{\bar{a}( \bar{x})}D_{\bar{x}} \bar{K}.\\
\end{array}
\end{equation}
Now, we need to answer two main questions. How can one obtain the
function $\bar{a}'(\bar{x})$ of the equivalent canonical form
$\bar{u}'''=\bar{a}^3( \bar{x})  \bar{u}$ by use of the obtained
invariant $K$? How can one determine the point transformation
(\ref{cccc}) that does the reduction of $u''' = f(x,u,u',u'')$
to the equivalent canonical form $\bar{u}'''=\bar{a}^3( \bar{x})
\bar{u}$ efficiently?

To answer these two questions, we need to look at the invariant
coframe (\ref{b36}) $\textrm{mod}~(\omega^1,\omega^2,\omega^3)$
and introduce some auxiliary functions. The symmetrical version of
the Cartan formulation $\textrm{mod}~(\omega^1,\omega^2,\omega^3)$
is
\begin{equation}\label{b37}
\begin{array}{l}
  J dx=\bar{J} d\bar{x}, \\
  \frac{d a_1}{a_1}-\frac{1}{3}\frac{s_4}{J}dx=\frac{d \bar{a}_1}{\bar{a}_1}-\frac{1}{3}\frac{\bar{s}_4}{\bar{J}}d\bar{x}.\\
\end{array}
\end{equation}
Inserting the point transformation (\ref{cccc}) into (\ref{b37})
and by use of  $\bar{J}=\bar{a}( \bar{x}),~\bar{s}_4=3D_{\bar{x}}
\bar{a}( \bar{x})$ for $\bar{f}= \bar{a}^3( \bar{x})  \bar{u}$ and
$\bar{a_1}=1$, results in
\begin{equation}\label{b39}
\begin{array}{l}
  J=\bar{a}( \bar{x}) D_x \phi, \\
  \frac{D_x a_1}{a_1}-\frac{1}{3}\left(\frac{3 D_x J-J f_q}{J}\right)=-\frac{D_{\bar{x}}\bar{a}( \bar{x})}{\bar{a}( \bar{x})}D_x \phi,\\
\end{array}
\end{equation}
where the auxiliary function $a_1(x,u,p)=\frac{\phi_x\psi_u   -
\phi_u \psi_x}{D_x \phi}$.

Invoking system (\ref{b39}), we introduce another auxiliary function
$H(x,u)$ as
\begin{equation}\label{z9}
\begin{array}{l}
     H=\frac{D_{\bar{x}}\bar{a}( \bar{x})}{\bar{a}( \bar{x})^2}=\frac{1}{3J}\left(\frac{3 D_x J-J f_q}{J}\right)-\frac{1}{J}\frac{D_x a_1}{a_1}.\\
\end{array}
\end{equation}
Thus, the first equation in the system (\ref{z8}) reads
\begin{equation}\label{z10}
\begin{array}{l}
K=\frac{2}{J}D_xH+H^2.\\
\end{array}
\end{equation}
By utilizing the fact $D_x=(D_x \phi)D_{\bar{x}}=\frac{J}{\bar{a}(
\bar{x})}D_{\bar{x}}$, equations (\ref{z9}) can be rewritten  as
\begin{equation}\label{z11}
\begin{array}{l}
D_x b =(J H) b,\\
D_x a_1 =\left(\frac{D_x J}{J}-\frac{1}{3}f_q-J H\right)a_1,\\
\end{array}
\end{equation}
where $\bar{a}( \bar{x})=\bar{a}(\phi(x,u))=b(x,u).$

Thus, one can start solving equation (\ref{z10}) for
$H(x,u)$. Then solve (\ref{z11}) for $b(x,u)$ and $a_1(x,u,p)$.
Finally, one can determine the point transformation (\ref{cccc})
by solution of the following  system
\begin{equation}\label{z12}
\begin{array}{l}
D_x \phi=\frac{J}{b}, \\
\phi_x\psi_u   - \phi_u \psi_x=\frac{J a_1}{b}.\\
\end{array}
\end{equation}
This proves the following theorem.
\begin{theorem}\label{th1}
The necessary and sufficient conditions for equivalence of a
third-order equation $ u''' = f(x,u,u',u'')$ via {\it point transformation} (\ref{cccc}) to the canonical form
$\bar{u}'''=\bar{a}^3( \bar{x})  \bar{u}$, with four point
symmetries,  are the
identically vanishing of the relative invariants
\begin{equation}\label{b40}
\begin{array}{l}
  I_1  = f_{q,q,q}  \hfill \\
  I_2  = f_{q,q} ^2  + 6\,f_{p,q,q}  \hfill \\
  I_4  = J_q \,\,\,\,\,\,\,\,\,\,\,\,\,\,\,\,\,\,\,\,\,\,\,\,\,\,\,\,\,\,\,\,\,\,\,\,\,\,\,\,\,\,\,\,\,\,\,\,\,\,\,\,\,\,\,\,\,\,\,\,\,\,\,\,\,\,\,\,\,\,\,\,\,\,\,\,\,\,\,\,\,\,\,\,\,\,\,\,\,\,\,\,\,\,\,\,\,\,\,\,\, \hfill \\
  I_5  = f_{q,q} \,J - 6\,J_p \,\,\,\,\,\,\,\,\,\,\,\,\,\,\,\,\,\,\,\,\,\,\,\,\,\,\,\,\,\,\, \hfill \\
  I_6  = J_u  - D_x J_p  \hfill \\
  I_7  = f_{q,q} \left( {9\,f_p  + f_q ^2  - 3\,D_x f_q } \right) - 9\,f_{p,p}  + 18\,f_{u,q}  - 6\,f_q f_{p,q}  \hfill \\
  I_9  =K_q\\
  I_{10}=K_p\\
  I_{11}=f_{q,q}D_xK-6K_u\\
\end{array}
\end{equation}
where
\begin{equation}\label{b41}
\begin{array}{l}
I_3  =J^3= \frac{1}{54} \left( {4\,f_q ^3  + 18\,f_q \left( {f_p  - D_x f_q } \right) + 9\,D_x^2 f_q    - 27\,D_x f_p+ 54\,f_u  }\right)  \ne 0,\\
I_8 = \frac{1}{3} \left( \left( {f_q ^2  + 3\,f_p  - 3\,D_x f_q \,} \right)\,J^2 \, + 6\,J\,D_x^2 J - 9\,\left( {D_x J} \right)^2\right)\ne 0, \\
K=\frac{I_8}{J^4},\\
D_xK\ne0,\\
\end{array}
\end{equation}
Given that the the system of relative invariants (\ref{b40}) is
zero and conditions (\ref{b41}) are satisfied, the linearizing
point transformation (\ref{cccc}) is defined by
\begin{equation}\label{b42}
\begin{array}{l}
D_x \phi=\frac{J}{b}, \\
\phi_x\psi_u   - \phi_u \psi_x=\frac{J a_1}{b},\\
\end{array}
\end{equation}
where  $H(x,u), b(x,u), a_1(x,u,p)$ are auxiliary functions given
by
\begin{equation}\label{b43}
\begin{array}{l}
\frac{2}{J}D_xH+H^2=K,\\
D_x b =(J H) b,\\
D_x a_1 =\left(\frac{D_x J}{J}-\frac{1}{3}f_q-J H\right)a_1.\\
\end{array}
\end{equation}
Finally, the function $\bar{a}( \bar{x})$ of the resulting
canonical form is given by the equation
\begin{equation}\label{b44}
\begin{array}{l}\bar{a}( \bar{x})=\bar{a}(\phi(x,u))=b(x,u).\\
\end{array}
\end{equation}
\end{theorem}
\section{Illustration of Theorem \ref{th1}}
\begin{example}\rm \label{ex1}
Consider the nonlinear ODE
\begin{equation}\label{n1}
\begin{array}{c}
u'''=3\frac{{u''}^2}{u'}-xu^3 {u'}^4.\\
\end{array}
\end{equation}
The relative invariants (\ref{b40}) are equal to zero while the relative
invariants (\ref{b41}) are not equal to zero for the function
\begin{equation}\label{n2}
\begin{array}{ll}
f(x,u,p,q)=3\frac{{q}^2}{p}-xu^3 {p}^4.\\
\end{array}
\end{equation}
As a consequence, this equation via  point transformation (\ref{cccc}) is equivalent to the canonical
form $\bar{u}'''=\bar{a}^3( \bar{x})  \bar{u}$, with four point
symmetries. We outline the
steps for constructing this point transformation  and finding the
function $\bar{a}( \bar{x})$ of the resulting canonical form.
\begin{description}
\item[Step 1] Solve the system (\ref{b43}) for the auxiliary functions
$H(x,u), b(x,u), a_1(x,u,p)$ as shown below.

By using the total derivative definition and comparing the
coefficients in the first equation of the system (\ref{b43})
\begin{equation}\label{n3}
\begin{array}{lll}
\frac{2}{up}D_xH+H^2=-\frac{3}{u^4},\\
\end{array}
\end{equation}
yields that $H=H(u)$ which satisfies the Riccati equation
\begin{equation}\label{n4}
\begin{array}{l}
2\frac{\partial H}{\partial u}+uH^2=-\frac{3}{u^3},\\
\end{array}
\end{equation}
with a solution  $H=\frac{1}{u^2}$.

The second equation of the system (\ref{b43}) reads
\begin{equation}\label{n5}
\begin{array}{l}
D_x b = (\frac{p}{u}) b.\\
\end{array}
\end{equation}
Now by utilizing the total derivative definition and comparing the
coefficients of equation (\ref{n5}), results in
\begin{equation}\label{n6}
\begin{array}{ll}
\frac{\partial b}{\partial u}=(\frac{1}{u}) b,&\frac{\partial b}{\partial x}=0,\\
\end{array}
\end{equation}
which has a solution  $b=u$.

Now, the third equation of the system (\ref{b43}) is
\begin{equation}\label{n7}
\begin{array}{lll}
D_x a_1 =-\left(\frac{q}{p}\right)a_1.\\
\end{array}
\end{equation}
Now, utilizing the total derivative definition and comparing the
coefficients of of the variable $q$ in equation (\ref{n7}), yield
\begin{equation}\label{n8}
\begin{array}{l}
\frac{\partial a_1}{\partial p}=-\frac{1}{p} a_1,\\
\end{array}
\end{equation}
which gives $a_1=\frac{F(x,u)}{p}$. Substituting $a_1$ back into
equation (\ref{n7}) and then comparing the coefficients, one arrives at
\begin{equation}\label{n9}
\begin{array}{ll}
\frac{\partial F}{\partial u}=0, &\frac{\partial F}{\partial x}=0,\\
\end{array}
\end{equation}
which gives a solution $F=1$. Thus, $a_1=\frac{1}{p}$ is a
solution of equation (\ref{n7}).

\item[Step 2] Obtain  the linearizing transformation (\ref{cccc}) by
 solving the system (\ref{b42}) as shown below.

The first equation of the system (\ref{b42}) reads
\begin{equation}\label{n11}
\begin{array}{l}
D_x\phi=p,\\
\end{array}
\end{equation}
which has a solution $\phi(x,u)=u$.

Now, the second equation of the system (\ref{b42}) is
\begin{equation}\label{n12}
\begin{array}{l}
\psi_x=-1,\\
\end{array}
\end{equation}
which has a solution $\psi(x,u)=-x$.

Finally, the function $\bar{a}( \bar{x})$ of the resulting
canonical form is given by the equation (\ref{b44}) as $\bar{a}(
\bar{x})=\bar{a}(u)=u$. Therefore, $\bar{a}( \bar{x})=\bar{x}$ and
the canonical form $\bar{u}'''=\bar{x}^3 \bar{u}$  can be found
for the ODE (\ref{n1}) via the point transformation
$$\bar{x}=u, \bar{u}=-x.$$
\end{description}
\end{example}
\section{Equivalence to to Yumaguzhin's canonical form}
Yumaguzhin \cite{Yum} showed that a  scalar  third-order linear
ODE with a four-dimensional symmetry algebra is equivalent via point
transformation to
\begin{equation}\label{e1}
u'''=\left(h{g}^2-2\left(\frac{g'}{g}\right)'+\left(\frac{g'}{g}\right)^2\right) u'+\frac{1}{2}\left(\left( hg^2-2\left(\frac{g'}{g}\right)'+\left(\frac{g'}{g}\right)^2 \right)'- g^3 \right)u,\\
\end{equation}
where $g$ is an arbitrary nowhere vanishing smooth function of $x$
and $h$ is  a function of $x$ with non-zero first derivative.

Using the point transformation $\bar{x}=\pm h, \bar{u}=\pm h' u$,
equation (\ref{e1}) is transformed to the canonical form
\begin{equation}\label{ee1}
\bar{u}'''=\left(\pm\bar{x}{\bar{g}}^2-2\left(\frac{\bar{g}'}{\bar{g}}\right)'+\left(\frac{\bar{g}'}{\bar{g}}\right)^2\right) \bar{u}'+\frac{1}{2}\left(\left(\pm \bar{x}{\bar{g}}^2-2\left(\frac{\bar{g}'}{\bar{g}}\right)'+\left(\frac{\bar{g}'}{\bar{g}}\right)^2 \right)'- \bar{g}^3 \right)\bar{u}.\\
\end{equation}
Now, if we scale $J$ to be
\begin{equation}\label{e0}
\begin{array}{l}
I_3  =J^3= -\frac{1}{27} \left( {4\,f_q ^3  + 18\,f_q \left( {f_p  - D_x f_q } \right) + 9\,D_x^2 f_q    - 27\,D_x f_p+ 54\,f_u  }\right),\\
\end{array}
\end{equation}
then $\bar{J}=\bar{g}( \bar{x}),~\bar{s}_4=3D_{\bar{x}} \bar{g}(
\bar{x})$ for $\bar{f}$ of the canonical form (\ref{ee1}).

Inserting the point transformation (\ref{cccc}) into (\ref{b37})
and using $\bar{J}=\bar{g}( \bar{x}),~\bar{s}_4=3D_{\bar{x}}
\bar{g}( \bar{x})$ and $\bar{a_1}=1$, gives rise to
\begin{equation}\label{e2}
\begin{array}{l}
  J=\bar{g}( \bar{x}) D_x \phi, \\
  \frac{D_x a_1}{a_1}-\frac{1}{3}\left(\frac{3 D_x J-J f_q}{J}\right)=-\frac{D_{\bar{x}}\bar{g}( \bar{x})}{\bar{g}( \bar{x})}D_x \phi,\\
\end{array}
\end{equation}
where the auxiliary function $a_1(x,u,p)=\frac{\phi_x\psi_u   -
\phi_u \psi_x}{D_x \phi}$.

Moreover, the invariants $K$ and $\frac{1}{\bar{J}}D_{x} K$ can be
given  simply as
\begin{equation}\label{e3}
\begin{array}{l}
K=\bar{K}=\pm\bar{x},\\
\frac{1}{\bar{J}}D_{x} K=\pm\frac{1}{\bar{g}( \bar{x})}.\\
\end{array}
\end{equation}

This proves another version of Theorem \ref{th1}.
\begin{theorem}\label{th2}
The necessary and sufficient conditions for equivalence of a
third-order equation $ u''' = f(x,u,u',u'')$ via {\it point
transformation} (\ref{cccc})  to the canonical form (\ref{ee1}),
with four point symmetries, are the identically vanishing of the
relative invariants
\begin{equation}\label{e4}
\begin{array}{l}
  I_1  = f_{q,q,q}  \hfill \\
  I_2  = f_{q,q} ^2  + 6\,f_{p,q,q}  \hfill \\
  I_4  = J_q \,\,\,\,\,\,\,\,\,\,\,\,\,\,\,\,\,\,\,\,\,\,\,\,\,\,\,\,\,\,\,\,\,\,\,\,\,\,\,\,\,\,\,\,\,\,\,\,\,\,\,\,\,\,\,\,\,\,\,\,\,\,\,\,\,\,\,\,\,\,\,\,\,\,\,\,\,\,\,\,\,\,\,\,\,\,\,\,\,\,\,\,\,\,\,\,\,\,\,\,\, \hfill \\
  I_5  = f_{q,q} \,J - 6\,J_p \,\,\,\,\,\,\,\,\,\,\,\,\,\,\,\,\,\,\,\,\,\,\,\,\,\,\,\,\,\,\, \hfill \\
  I_6  = J_u  - D_x J_p  \hfill \\
  I_7  = f_{q,q} \left( {9\,f_p  + f_q ^2  - 3\,D_x f_q } \right) - 9\,f_{p,p}  + 18\,f_{u,q}  - 6\,f_q f_{p,q}  \hfill \\
  I_9  =K_q\\
  I_{10}=K_p\\
  I_{11}=f_{q,q}D_xK-6K_u\\
\end{array}
\end{equation}
where
\begin{equation}\label{e5}
\begin{array}{l}
I_3  =J^3= -\frac{1}{27} \left( {4\,f_q ^3  + 18\,f_q \left( {f_p  - D_x f_q } \right) + 9\,D_x^2 f_q    - 27\,D_x f_p+ 54\,f_u  }\right)  \ne 0,\\
I_8 = \frac{1}{3} \left( \left( {f_q ^2  + 3\,f_p  - 3\,D_x f_q \,} \right)\,J^2 \, + 6\,J\,D_x^2 J - 9\,\left( {D_x J} \right)^2\right)\ne 0, \\
K=\frac{I_8}{J^4}.\\
D_xK\ne0,\\
\end{array}
\end{equation}
Given that the the system of relative invariants (\ref{e4}) is
zero and conditions (\ref{e5}) are satisfied, the linearizing
point transformation (\ref{cccc}) is defined by $\phi=\pm K$ and
the linear equation
\begin{equation}\label{e6}
\begin{array}{l}
K_x\psi_u   - K_u \psi_x=a_1 D_x K,\\
\end{array}
\end{equation}
where  $a_1(x,u,p)$ are auxiliary functions given by
\begin{equation}\label{e7}
\begin{array}{l}
D_x a_1 =\left(\frac{D_x J}{J}-\frac{1}{3}f_q- \frac{D_{x} K}{J}  D_x \left(\frac{J}{D_{x} K}\right)\right)a_1.\\
\end{array}
\end{equation}
Finally, the function $\bar{g}( \bar{x})$ of the resulting
canonical form is given by the equation
\begin{equation}\label{e8}
\begin{array}{l}\bar{g}( \bar{x})=\bar{g}(\pm K)=\pm \frac{J}{D_{x} K}.\\
\end{array}
\end{equation}
\end{theorem}
\section{Illustration of Theorem \ref{th2}}
\begin{example}\rm \label{ex2}
Consider the nonlinear ODE
\begin{equation}\label{m1}
\begin{array}{c}
u'''=3\frac{{u''}^2}{u'}-xu^3 {u'}^4.\\
\end{array}
\end{equation}
As explained in example \ref{ex1}, this equation admits four point
symmetries. Using (\ref{e8}), since $J=- \sqrt[3]{2}~u p$,
$I_8=-3\sqrt[3]{4}~ p^4$, $K=-\frac{3}{\sqrt[3]{4}~ }u^{-4}$,
$\frac{J}{D_{x} K}= -\frac{1}{6}u^6$, the ODE (\ref{m1}) is
equivalent via point transformation (\ref{cccc}) to Yumaguzhin's
canonical form (\ref{ee1}) with positive sign and the complex
function $\bar{g}( \bar{x})=i
\frac{\sqrt{3}}{4}~\bar{x}^{-\frac{3}{2}}$. Notwithstanding, the ODE
(\ref{m1}) is equivalent via point transformation (\ref{cccc}) to
Yumaguzhin's canonical form (\ref{ee1}) with negative sign and
real function $\bar{g}(
\bar{x})=\frac{\sqrt{3}}{4}~\bar{x}^{-\frac{3}{2}}$.

We outline the steps for constructing the point transformation for
Yumaguzhin's canonical form (\ref{ee1}) with negative sign and the
real function $\bar{g}(
\bar{x})=\frac{\sqrt{3}}{4}~\bar{x}^{-\frac{3}{2}}$.
\begin{description}
\item[Step 1] Solve the equation (\ref{e7}) for the auxiliary functions
$a_1(x,u,p)$.

The equation (\ref{e7}) reads
\begin{equation}\label{m2}
\begin{array}{lll}
D_x a_1 =-\left(\frac{q}{p}+5\frac{p}{u}\right)a_1.\\
\end{array}
\end{equation}
Now, utilizing the total derivative definition and comparing the
coefficients of of the variable $q$ in equation (\ref{m2}), yield
\begin{equation}\label{m3}
\begin{array}{l}
\frac{\partial a_1}{\partial p}=-\frac{1}{p} a_1,\\
\end{array}
\end{equation}
which gives $a_1=\frac{F(x,u)}{p}$. Substituting $a_1$ back into
equation (\ref{m2}) and then comparing the coefficients, result in
\begin{equation}\label{m4}
\begin{array}{ll}
\frac{\partial F}{\partial u}=-\frac{5}{u}F, &\frac{\partial F}{\partial x}=0,\\
\end{array}
\end{equation}
which gives rise to a solution $F=\frac{1}{u^{5}}$. Thus,
$a_1=\frac{1}{u^{5}p}$ is a solution of equation (\ref{m2}).

\item[Step 2]Solve the system (\ref{e6}) for $\psi(x,u)$.

The equation (\ref{e6}) reads
\begin{equation}\label{m5}
\begin{array}{lll}
\psi_x=-\frac{1}{u^{5}} ,\\
\end{array}
\end{equation}
which has a solution $\psi(x,u)=-\frac{x}{u^{5}}$.
\end{description}
Therefore, the Yumaguzhin's canonical form (\ref{ee1}) with
negative sign and the real function $\bar{g}(
\bar{x})=\frac{\sqrt{3}}{4}~\bar{x}^{-\frac{3}{2}}$  can be
obtained for the ODE (\ref{m1}) via the point transformation
\begin{equation}\label{m6}
\bar{x}=\frac{3}{\sqrt[3]{4}~ u^4}, \bar{u}=-\frac{x}{u^{5}}.\\
\end{equation}
\end{example}
\begin{remark}\label{r2}
The Yumaguzhin's canonical form (\ref{ee1}) with negative sign and
the real function
$\bar{g}(\bar{x})=\frac{\sqrt{3}}{4}~\bar{x}^{-\frac{3}{2}}$ can
be further transformed via the point transformation
\begin{equation}\label{m7}
\bar{x}=\frac{3}{\sqrt[3]{4}~ t^4}, \bar{u}=-\frac{y(t)}{t^{5}},\\
\end{equation} to the ODE
\begin{equation}\label{m8}
y'''(t)=t^3 y(t).\\
\end{equation}

It is clear that the composition of the transformation (\ref{m6})
and the transformation (\ref{m7}) is the transformation
\begin{equation}\label{m9}
x=y(t), u=t.\\
\end{equation}
\end{remark}
\section{Application for point symmetry group classification of third-order linear ODE}
\begin{corollary}\label{CL6.1}
The point symmetry group classification of the third-order linear
ODE
\begin{equation}\label{lkj}
\begin{array}{lll}
 u'''= a_1(x)u''+a_2(x)u'+a_3(x)u+a_4(x)\\
\end{array}
\end{equation}
is given as follows:
\begin{itemize}
    \item Seven point symmetries: $J=0$
    \item Five point symmetries: $J\ne0, D_x K=0$
    \item Four point symmetries: $J\ne0, D_x K\ne0$
\end{itemize}
where
\begin{equation}\label{b41}
\begin{array}{l}
f=a_1(x)q+a_2(x)p+a_3(x)u+a_4(x),\\
J^3= \frac{1}{54} \left( {4\,f_q ^3  + 18\,f_q \left( {f_p  - D_x f_q } \right) + 9\,D_x^2 f_q    - 27\,D_x f_p+ 54\,f_u  }\right),\\
I_8 = \frac{1}{3} \left( \left( {f_q ^2  + 3\,f_p  - 3\,D_x f_q \,} \right)\,J^2 \, + 6\,J\,D_x^2 J - 9\,\left( {D_x J} \right)^2\right), \\
K=\frac{I_8}{J^4}.\\
\end{array}
\end{equation}
\end{corollary}
\proof
It is well-known that the linear third-order ODEs admit
seven, five or four point symmetries.

By using Theorem \ref{TH1.1}, the ODE (\ref{lkj}) admits seven
point symmetries if and only if $J=0$. Also by utilizing Theorem
\ref{TH1.2}, the ODE (\ref{lkj}) admits five point symmetries if
and only if $J\ne0, D_x K=0$. Finally, by utilizing Theorem
\ref{th1}, the ODE (\ref{lkj}) admits four point symmetries if and
only if $J\ne0, D_x K\ne0$.
\endproof
\begin{example}\cite[page 285]{Michal2012}  \textbf{Regulation of a steam turbine}\\
The motion $x(t)$ of a steam supply control slide valve is
governed by the third-order differential equation
\begin{equation}\label{v1}
\begin{array}{lll}
m x'''+f x''+k x'+h \frac{\alpha}{I}x=0,\\
\end{array}
\end{equation}
where $m$ is the mass of the valve, $f$ is friction, $k$ is a
constant characterizing the properties of the slide valve spring,
$h$ is a constant depending on the dimensions of the equipment,
$\alpha$ is a proportionality constant relating the motion and the
acceleration of the control valve, and $I$ is the moment of
inertia of the turbine.

By using Corollary \ref{CL6.1}, the ODE (\ref{v1}) admits seven or
five point symmetries. There are seven point symmetries for
$\alpha=\frac{1}{27}{\frac {f{\it I}\, \left( 9\,km-2\,{f}^{2}
\right) }{{m}^{2}h}}$ and five point symmetries for $\alpha \ne
\frac{1}{27}{\frac {f{\it I}\, \left( 9\,km-2\,{f}^{2} \right)
}{{m}^{2}h}}$.
\end{example}

\begin{example}\cite[page 248]{Michal2012}  \textbf{Deflection of a Curved Beam having a constant or varying cross-section}\\
Consider the curved beam in the form of a circular arc with
diameter $a$ ($0\leq \theta \leq \gamma$, $\gamma>0$). In the case
of an equally distributed load on the beam, the corresponding
equation for the bending moment $M$ takes the form
\begin{equation}\label{v2}
\begin{array}{lll}
\frac{1}{a^3}\left(\frac{d^2}{d \xi^2}+1 \right)\frac{d M}{d \xi}+p \frac{d }{d \xi}\left(\frac{M}{B} \right)=0,\\
\end{array}
\end{equation}
where $p$ is the force acting on a unit length, $B$ is the
flexured rigidity (if the cross-section is constant, then B is
constant) and $\xi$ is an introduced independent variable given
explicitly in terms of the angular coordinate $\theta$ as
$\xi=\theta-\frac{\gamma}{2}$.

By using Corollary \ref{CL6.1}, the ODE (\ref{v2}) admits seven
point symmetries for constant flexured rigidity $B$ whereas it
admits five point symmetries for non-constant $B$ given by the
following differential constraint
\begin{equation}\label{v3}
\begin{array}{lll}
18\,{B_{{}}}^{3}{B_{{\xi}}}^{2}B_{{\xi,\xi}}-36\,{B_{{}}}^{2}{B_{{\xi}}}^{4}+18p{a}^{3}\,{B_{{}}}^{2}{B_{{\xi}}}^{2}B_{{\xi,\xi}}-9p{a}^{3}\,B_{{}}{B_{{\xi}}}^{4}+24\,B_{{}}{B_{{\xi}}}^{4}B_{{\xi,\xi}}\\
-12\,{B_{{}}}^{2}{B_{{\xi}}}^{2}{B_{{\xi,\xi}}}^{2}-72\,{B_{{}}}^{3}B_{{\xi}}B_{{\xi,\xi}}B_{{\xi,\xi,\xi}}+18\,{B_{{}}}^{3}{B_{{\xi}}}^{2}B_{{\xi,\xi,\xi,\xi}}-16\,{B_{{\xi}}}^{6}+56{B_{{}}}^{3}\,{B_{{\xi,\xi}}}^{3}=0.\\
\end{array}
\end{equation}
Finally, the ODE (\ref{v2}) admits four point symmetries for any
$B$ not satisfying the differential constraint (\ref{v3}).
\end{example}
\begin{remark}\label{r3}
One can verify that $B \left( \xi \right) =-p{a}^{3}{\frac
{{\xi}^{2}}{{\xi}^{2}-1}}$ is a solution for the differential
constraint (\ref{v3}). Therefore, the ODE (\ref{v2}) with this $B$
admits five point symmetries. Hence, it can be transformed to an ODE
with constant coefficients and its solution can be given
explicitly.
\end{remark}
\section{Conclusion}
In this work we have completed the solution to the invariant
characterization of the linearizable scalar third-order ODE
$u'''=f(x,u,u',u'')$ which admits a four-dimensional point
symmetry Lie algebra via the Cartan equivalence method. Earlier
\cite{DweikTahirFazal2,DweikTahirFazal1} a study was made of the
linearization of such ODEs for five- and seven-dimensional
symmetry algebras. Our invariant characterization here is given
compactly in terms of a function. Moreover, we provide a method
for the determination of the maps to the canonical form.
Furthermore, we gave relevant applications of our main results.
\subsection*{Acknowledgments}
Ahmad Y. Al-Dweik  is thankful to the King Fahd University of
Petroleum and Minerals for its continued support and excellent research
facilities. FMM  is grateful to the NRF of South Africa for research support.


\begin{thebibliography}{99}
\bibitem{DweikTahirFazal2}Ahmad Y. Al-Dweik, M. T. Mustafa, F. M. Mahomed, Invariant characterization of scalar third-order ODEs that admit the maximal point symmetry Lie algebra,     arXiv:1712.02387 [math.CA].
\bibitem{DweikTahirFazal1}Ahmad Y. Al-Dweik, M. T. Mustafa, F. M. Mahomed, Invariant characterization of third-order ODEs $u'''=f(x,u,u',u'')$ that admit a five-dimensional point symmetry Lie algebra, arXiv:1711.08138 [math.CA].
\bibitem{mah4} Mahomed, F. M. and Leach P. G. L., Symmetry Lie Algebras of $n$th Order Ordinary Differential Equations. {\it J Math Anal Applic}  {\bf151}, (1990), 80.
\bibitem{che} Chern, S.S., The geometry of the differential equation $y'''=F(x,y,y,y'')$, Sci. Rep. Nat. Tsing Hua Univ. 4 (1940), 97-111.
\bibitem{neu} Neut, S. and Petitot, M., La g\'eom\'etrie de l'\'equation $y'''=f(x,y,y',y'')$ {\it C.R. Acad. Sci. Paris S\'er I} {\bf 335}, (2002), 515-518.
\bibitem{gre} Grebot, G., The characterization of third order ordinary differential equations admitting a transitive fibre-preserving point symmetry group, {\it J. Math. Anal. Applic.} {\bf 206}, (1997), 364-388.
\bibitem{ibr} N.H. Ibragimov, S.V. Meleshko, Linearization of third-order ordinary differential equations by point and contact transformations, {\it J. Math. Anal. Appl.} {\bf308}, (2005), 266–-289.
\bibitem{Mahomed1996} N.H. Ibragimov and F.M. Mahomed, Ordinary differential equations, CRC Handbook of Lie Group Analysis of Differential Equations, vol. 3. N.H. Ibragimov ed., CRC Press, Boca Raton (1996) 191.
\bibitem{Mahomed2007}Mahomed F M 2007, Symmetry group classification of ordinary differential equations: Survey of some results, {\it Math Methods in the Applied Sciences}, {\bf 30}, 1995-2012.
\bibitem{Olver1995} Olver, P.J., Equivalence, Invariants and Symmetry, Cambridge University Press, Cambridge, 1995.
\bibitem{Yum}Yumaguzhin, Valeriy A. "Classification of 3rd order linear ODE up to equivalence." Differential Geometry and its Applications 6.4(1996): 343-350.
\bibitem{Michal2012}Gregus, Michal. Third order linear differential equations. Vol. 22. Springer Science \& Business Media, 2012.
\end{thebibliography}
\end{document}